\documentclass[12pt]{article}

\usepackage{amsmath,amssymb,amsfonts,theorem,makeidx,latexsym,epsfig,,subfigure}

\newtheorem{defn}{Definition}[section]

\newtheorem{lemma}[defn]{Lemma}

{\theorembodyfont{\rmfamily}

\newtheorem{ex}[defn]{Example}}

\newtheorem{thm}[defn]{Theorem}

\newtheorem{prop}[defn]{Proposition}

\newtheorem{cor}[defn]{Corollary}

\newtheorem{rem}[defn]{Remark}

\numberwithin{equation}{section}

\newcommand{\ltr}{ L^2(\mathbb R) }

\newcommand{\mn}{\mathbb N}

\newcommand{\mr}{\mathbb R}

\newcommand{\mz}{\mathbb Z}

\newcommand{\mc}{\mathbb C}

\def\bp{{\noindent\bf Proof. \ }}

\def\ep{\hfill$\square$\par\bigskip}

\def\bqs{\begin{equation}}

\def\eqs{\tag*{$\square$}\end{equation}\par\bigskip}

\def\la{\langle}

\def\ra{\rangle}

\def\ga{\gamma}

\def\supp{\text{supp}}

\def\bop{\begin{op}\rm}

\def\eop{\end{op}}

\def\bee{\begin{eqnarray}}

\def\ene{\end{eqnarray}}

\def\bes{\begin{eqnarray*}}

\def\ens{\end{eqnarray*}}

\def\bei{\begin{itemize}}

\def\eni{\end{itemize}}

\def\bt{\begin{thm}}

\def\et{\end{thm}}

\def\bc{\begin{cor}}

\def\ec{\end{cor}}

\def\bpr{\begin{prop}}

\def\epr{\end{prop}}

\def\bl{\begin{lemma}}

\def\el{\end{lemma}}

\def\bd{\begin{defn}}

\def\ed{\end{defn}}

\def\bex{\begin{ex}}

\def\enx{\end{ex}}

\def\bfi{\begin{fig}}

\def\efi{\end{fig}}

\def\inr{\int_{-\infty}^\infty}

\def\inrd{\int_{\mrd}}

\def\mzd{{\mathbb Z}^d}

\def\mrd{{\mathbb R}^d}

\newcommand{\ltrd}{ L^2({\mathbb R}^d) }

\def\rt{R(\theta)}
\def\rtp{R(\theta^\prime)}

\def\cam{{\cal A}}

\DeclareMathOperator*{\esssup}{ess\,sup}
\DeclareMathOperator*{\essinf}{ess\,inf}

\title{Construction of scaling partitions of unity }

\date{}

\author{Ole Christensen and Say Song Goh}

\begin{document}

\maketitle

\begin{abstract} Partitions of unity in $\mrd$ formed by (matrix) scales of a fixed function appear in many
parts of harmonic analysis, e.g., wavelet analysis and the analysis
of Triebel-Lizorkin spaces. We give a
simple characterization of the functions and matrices yielding such a partition of unity.
For invertible expanding matrices, the characterization leads to easy ways of constructing appropriate functions
with attractive properties like high regularity and small support.
We also discuss a class of integral transforms that map functions
having the  partition of unity property to functions with the same property. The one-dimensional version of the transform allows a direct definition of a class of
nonuniform splines with properties that are parallel to those of the classical
B-splines. The results are illustrated with the construction of dual pairs of wavelet frames.
%

\end{abstract}

\section{Introduction} A function $g: \mrd \to \mc$ is said to have the {\it $($scaling$)$ partition
of unity property} with respect to a real invertible $d\times d$ matrix $\cam$ if
\bee \label{60106cf} \sum_{j=-\infty}^\infty g(\cam^j\ga)=1, \forall \ga \in \mrd \setminus \{0\}.\ene Partitions of unity of this form appear in several
parts of analysis, e.g., wavelet analysis and the theory for
Triebel-Lizorkin spaces,
and the question of how to construct them has attracted some attention. In particular,
this issue comes up in connection with the analysis of tight wavelet frames in $\ltrd$
\cite{DDG3} and the more general case of dual wavelet frame pairs \cite{Lem1,Lem2}.

In this paper we will give a surprisingly simple characterization of the scaling partition of
unity property. In the case where $\cam$ is an expanding matrix, 
i.e., a real
matrix $\cam$ such that all its eigenvalues have absolute value strictly greater
than one, the characterization leads to easy ways of constructing appropriate functions $g$
with attractive properties like high regularity and small support.
Under certain
conditions, nonnegativity of the function $g$ can be guaranteed.
We also discuss a class of integral transforms that can be used to generate functions
with the  partition of unity property. The one-dimensional version of the transform
leads in a natural fashion to a definition of a recursively
given family of
nonuniform splines. These splines have some similarities with the classical B-splines:
their regularity and support grow with the order, and they satisfy the de Boor
recursion formula.  However, there are also differences: all the splines have support
within $[-1,1],$ and they satisfy a scaling partition of unity condition
instead of the translation partition of unity condition.
Finally, the key results  are applied to the construction of dual pairs of wavelet frames.

The paper is organized as follows.  In Section \ref{60129a}, we characterize the scaling partition
of unity condition and provide explicit and easily verifiable sufficient conditions
in the case where $\cam$ is an expanding matrix. 
Section \ref{sec3} deals with
the above mentioned one-dimensional integral transform and its lifting to higher dimensions. Finally,
Section \ref{60129c} applies the results to obtain easy constructions of wavelet frames
in $\ltrd$ and their associated dual frames.


\section{Characterization of the  partition of unity property} \label{60129a}

We first establish a characterization of the scaling partition of unity property. Despite
its simplicity we have not been able to find it stated in the literature.

\bt \label{60122a}
Consider a function $g: \mrd \to \mc$ and any real $d\times d$ matrix $\cam.$ Then the following hold:

\bei \item[{\rm (i)}] Assume that the infinite
series $\sum_{j=-\infty}^\infty g(\cam^j\ga)$ is convergent for all $\ga \in \mrd \setminus \{0\}.$  Then there is a function $\varphi: \mrd \to \mc$ such that
\bee \label{60106bf} g(\ga)= \varphi(\ga)- \varphi(\cam\ga), \, \forall \ga \in \mrd \setminus \{0\}.\ene
\item[{\rm (ii)}] On the other hand, take any function $\varphi: \mrd \to \mc$ such that
\eqref{60106bf} holds. Then, fixing any
 $\ga\in \mrd \setminus \{0\}$, the series $\sum_{j=-\infty}^\infty g(\cam^j\ga)$ is convergent
if and only if the two limits  $\lim_{N \to \pm \infty} \varphi(\cam^N\ga)$ exist.
\item[{\rm (iii)}] Take again any function $\varphi: \mrd \to \mc$ such that
\eqref{60106bf} holds. Then the partition of unity condition \eqref{60106cf}  holds if and only if
the two limits  $\lim_{N \to \pm \infty} \varphi(\cam^N\ga)$ exist and
\bes 
     \lim_{N \to - \infty} \varphi(\cam^N\ga)- \lim_{N \to  \infty} \varphi(\cam^N\ga)  =1 \ens
     for all $\ga \in \mrd \setminus \{0\}.$
\eni
\et

\bp For the proof of (i), assume that the infinite
series $\sum_{j=-\infty}^\infty g(\cam^j\ga)$ is convergent for all $\ga \in \mrd \setminus \{0\}.$ Then
\bes g(\ga) =  \sum_{j=0}^\infty g(\cam^j\ga)-\sum_{j=1}^\infty g(\cam^j\ga) =
\sum_{j=0}^\infty g(\cam^j\ga)-\sum_{j=0}^\infty g(\cam^j\cam\ga).
\ens Taking now
$\varphi(\ga):= \sum_{j=0}^\infty g(\cam^j\ga), \, \ga \in \mrd \setminus \{0\},$
yields the result.
For the proof of (ii), by direct calculation and for any $M,N\in \mn,$
\bes \sum_{j=-M}^N g(\cam^j\ga)
& = & [\varphi(\cam^{-M}\ga)-\varphi(\cam^{-M+1}\ga)] + [\varphi(\cam^{-M+1}\ga)-\varphi(\cam^{-M+2}\ga)]   \\[-1ex] & \ & +
\cdots + [\varphi(\cam^{N}\ga)-\varphi(\cam^{N+1}\ga)]\\ & = & \varphi(\cam^{-M}\ga)-\varphi(\cam^{N+1}\ga).\ens
Then (ii) follows immediately; and (iii) is a consequence of (ii).
\ep

Note that the function $\varphi$ satisfying \eqref{60106bf} for a given function $g$
is not unique. In the sequel $\varphi$ will denote {\it any} such function, not necessarily
the one constructed in the proof of Theorem \ref{60122a}.

Via Theorem \ref{60122a}, we can now show that any expanding matrix $\cam$ leads to the
partition of unity property for a large class of functions $g.$  The following result and its proof hold whenever $|| \cdot ||$ denotes an arbitrary norm
on $\mrd.$

\bpr \label{60108f} Let $\cam$ be any invertible expanding $d\times d$ matrix, and consider any function $\varphi: \mrd \to \mc$
which is continuous at $\ga=0$ and satisfies the conditions that $\varphi(0)=1$ and
$\lim_{||\ga||\to \infty} \varphi(\ga)=0.$ Then the function
$g(\ga):= \varphi(\ga)- \varphi(\cam\ga)$ satisfies the partition of unity condition
\eqref{60106cf}. \epr

\bp By Lemma 5.2 in \cite{HLW}, an invertible matrix $\cam$ is expanding if and only if there
exist constants $C\in (0,1]$ and $\alpha>1$ such that
\bee \label{60113a}   || \cam^N \ga || \ge C \alpha^N ||\ga||\ene for all $\ga\in \mrd$ and
$N\in \mn \cup \{0\}.$ Thus, the assumption $\lim_{||\ga||\to \infty} \varphi(\ga)=0$ immediately implies that
$\lim_{N \to  \infty} \varphi(\cam^N\ga)=0$ for all $\ga \in \mrd \setminus\{0\}.$
Replacing $\ga$ by
$\cam^{-N} \ga$ in the inequality  \eqref{60113a}
shows that $||\cam^{-N}\ga|| \le C^{-1} \alpha^{-N}||\ga||$
for all $\ga\in \mrd$ and $N\in \mn \cup \{0\};$ thus the assumptions imply that
$\lim_{N \to - \infty} \varphi(\cam^N\ga)=1.$ The result now follows from Theorem \ref{60122a}. \ep

\bex \label{60106g} We first give an example of a partition of unity based on a
diagonal matrix, and then a construction that works for arbitrary expanding matrices.

\vspace{.1in}\noindent (i) Consider an even, continuous and nonnegative function $k:\mr \to  \mr$
such that $\int_0^\infty k(t)\,dt=1.$ Then the function
\bes \varphi(\ga):= \int_\ga^\infty k(t)\, dt, \, \ga \in \mr,\ens satisfies the conditions in Proposition \ref{60108f}. Thus, for any $a>1,$ the function
\bes g(\ga)= \varphi(\ga)- \varphi(a\ga)= \int_\ga^{a\ga} k(t)\, dt\ens
satisfies the partition of unity condition $\sum_{j\in \mz} g(a^j\ga)=1, \,
\ga \in \mr \setminus \{0\}.$ Clearly, $g\in C^1(\mr).$ Note that
for any choice of a norm $|| \cdot ||$ on $\mrd,$ the
function $g$ can be lifted to a radial function $\widetilde{g}: \mrd \to \mr,$
by defining $\widetilde{g}(\ga):= g(||\ga||), \, \ga \in \mrd;$ the function $\widetilde{g}$ satisfies the
partition of unity condition with respect to the $d \times d$ diagonal matrix $\cam=aI.$

\vspace{.1in}\noindent(ii) Let $|| \cdot||$  be the Euclidean norm on $\mrd.$
The function $\varphi(\ga):=  e^{-||\ga||^2}, \, \ga \in \mrd,$ satisfies the conditions in Proposition \ref{60108f}. Thus, for any expanding $d\times d$ matrix $\cam,$ the function
\bes g(\ga)= \varphi(\ga)- \varphi(\cam \ga)=  e^{-||\ga||^2} - e^{-||\cam\ga||^2}    \ens
satisfies the partition of unity condition \eqref{60106cf}. Clearly, $g\in C^\infty(\mrd).$
\ep \enx

Proposition \ref{60108f} makes it easy to construct partitions of unity for arbitrary
expanding matrices $\cam.$ Furthermore, several properties of the generating function $g$ can
be controlled directly in terms of the function $\varphi,$ e.g., regularity and support.
We now prove that nonnegativity of $g$ can also be guaranteed by choosing
$\varphi$ to be a radial function with respect to a given norm $||\cdot ||$
on $\mrd$:

\bpr \label{60113d} Let $||\cdot ||$ be an arbitrary norm on $\mrd$ and consider an expanding $d \times d$ matrix $\cam$ such that $||\ga||\le ||\cam\ga||$ for all $\ga \in \mrd.$
Let $r: [0,\infty) \to \mr$ denote a continuous decreasing function such that
$r(0)=1$ and $r(s) \to 0$ as $s\to \infty.$ Letting $\varphi(\ga):=r(||\ga||), \, \ga \in \mrd,$ the
function $g(\ga)=\varphi(\ga)-\varphi(\cam\ga)$ has the following properties:
\bei
\item[{\rm (i)}] $g\ge0.$
\item[{\rm (ii)}] $\sum_{j\in \mz} g(\cam^j \ga)=1, \forall \ga \in \mrd \setminus \{0\}.$
\item[{\rm (iii)}] There exists a constant $C>0$ such that
\bee \label{60113c} C \le\sum_{j\in \mz} |g(\cam^j \ga)|^2 \le 1, \, \forall \ga \in \mrd \setminus \{0\}.\ene \eni \epr

\bp Since the function $r$ is decreasing, (i) follows immediately
from the assumption that $||\ga||\le ||\cam\ga||$ for all $\ga \in \mrd.$ The
partition of unity (ii) follows from Proposition \ref{60108f}, so we only
need to prove (iii). In order to do so, the nonnegativity of $g$ and (ii) imply that
$0 \le g(\cam^j\ga)\le 1$ for every $j \in \mz$ and all $\ga \in \mrd \setminus \{0\};$ thus,
$\sum_{j\in \mz} |g(\cam^j \ga)|^2 \le \sum_{j\in \mz} g(\cam^j \ga)  =1.$

In order to prove the lower bound in \eqref{60113c}, let $\eta \in \mrd \setminus \{0\}.$
Then by (ii), there exists $j_\eta\in \mz$ such that $\epsilon_\eta:=g(\cam^{j_\eta}\eta)>0.$ Thus, we can
choose an open set $I_\eta$ containing $\eta$ such that $g(\cam^{j_\eta} \ga) \ge {\epsilon_\eta}/2$
for all $\ga \in I_\eta.$ Letting $B(0,1)$ denote the closed unit ball in $\mrd$
with respect to the norm $|| \cdot ||$,
the open sets $I_\eta, \eta \in B(0,1),$ form a cover of $B(0,1);$
thus, we can select a finite subcover, i.e.,
$B(0,1) \subset I_{\eta_1} \cup I_{\eta_2} \cup \cdots \cup I_{\eta_n}$
for some $\eta_1, \ldots, \eta_n \in B(0,1).$
It follows that for any $\ga \in \mrd$ with $|| \ga|| \le 1,$
$\ga$ must lie in $I_{\eta_\ell}$ for some $\ell \in \{1, \ldots, n \};$ thus,
\bes  \sum_{j\in \mz} |g(\cam^j \ga)|^2 \ge  |g(\cam^{j_{\eta_\ell}} \ga)|^2
\ge \frac14 \epsilon_{\eta_\ell}^2 \ge \frac14 \min \{ \epsilon_{\eta_1}^2, \ldots, \epsilon_{\eta_n}^2\}.\ens
This proves the lower bound in \eqref{60113c} for $\ga$ belonging to the closed unit ball in $\mrd.$ Taking
now an arbitrary $\ga \in \mrd \setminus \{0\},$ the argument in the proof of
Proposition \ref{60108f} shows that there exists $N\in \mn$ such that
$||\cam^{-N} \ga|| \le 1;$ thus, by a change of variable,
\bes  \sum_{j\in \mz} |g(\cam^j \ga)|^2= \sum_{j\in \mz} |g(\cam^{j+N}( \cam^{-N}\ga))|^2
=\sum_{j\in \mz} |g(\cam^{j}( \cam^{-N}\ga))|^2\ge \frac14 \min \{ \epsilon_{\eta_1}^2, \ldots, \epsilon_{\eta_n}^2\}.\ens
This completes the proof. \ep

The condition $||\ga||\le ||\cam\ga||, \ga \in \mrd,$ is clearly necessary for
the nonnegativity of $g(\ga)= \varphi(\ga)-\varphi(\cam\ga)$ whenever $\varphi$ is a function
of the type considered in Proposition \ref{60113d}. Note that the condition does not follow from $\cam$ being expanding,
as we shall see in the example below.

\bex \label{60113g} Take $|| \cdot||$ to be the Euclidean norm on $\mr^2$ and let
$\cam=
\begin{pmatrix} 0 & 2 \\ 3/4 & 0 \end{pmatrix}.$ The eigenvalues are $\pm \sqrt{3/2},$ so $\cam$ is
indeed expanding. However $\cam\begin{pmatrix} 1 \\ 0 \end{pmatrix} = \begin{pmatrix} 0 \\ 3/4 \end{pmatrix},$
so the condition $||\ga|| \leq ||\cam\ga||$ is clearly violated.
\ep \enx

\section{An integral transform preserving partitions of unity}
\label{sec3}

In this section, we consider certain integral
transforms that map a function $g$ having the scaling partition of unity property to another
function with the same property. We first discuss the transform on
$\mrd$ and then specialize to the one-dimensional case, where explicit calculations
are much easier. It turns out that the one-dimensional case leads to a
definition of a class of splines in a natural way.

\subsection{The integral transform on $\mrd$}
Fix a measurable function $g: \mrd \to \mc$ and
consider formally the integral operator $K_{g}$ that maps a function $f: \mrd \to \mc$ to
\bee \label{51224bda}  h(\ga)= (K_{g}f)(\ga):=\inrd f(t) g \big(\frac{\ga}{||t||}\big)\,dt, \, \ga \in \mrd,\ene
where $|| \cdot ||$ is an arbitrary norm on $\mrd.$
The set of functions $f$ for which the transform is well-defined clearly depends on the choice of the function $g.$  Typically, we assume that
$g$ is supported on
an annulus \bes 
a(R_1,R_2):= \{t\in \mrd\, \big| \, R_1 \le ||t|| \le R_2\}\ens for some $R_2>R_1>0$.
For example, if a function $f\in L^1(\mr)$  has support in an annulus
$a(R_1,R_2)$ and $g$ is a bounded function with support in an annulus
$a(R_3,R_4),$
then $h$ is well-defined and supported on the
annulus
$ a(R_3R_1,R_4R_2).$

The following proposition describes a case where the integral transform is well-defined
for all $f\in L^1(\mr)$ and generates a family of partitions of unity.

\bpr \label{51224eba} Let $g: \mrd \to \mc,$ and consider a real invertible $d\times d$ matrix $\cam$
such that
\bes 
\sum_{j\in \mz} g(\cam^j \ga)=1, \, \forall \ga \in \mrd \setminus \{0\},\ens
and there exists a constant $C>0$ for which
\bes \sum_{j\in \mz} |g(\cam^j \ga)| \le C, \, \forall \ga \in \mrd \setminus \{0\}.\ens
Then the integral transform $K_{g}$
in \eqref{51224bda} is well-defined for every $f\in L^1(\mrd),$ and
\bee \label{51224gba} \sum_{j\in \mz} h(\cam^j \ga)= \inrd f(t)\, dt, \, \forall \ga \in \mrd \setminus \{0\}.\ene 
In particular, if $f\in L^1(\mrd)$ is chosen such that $\int_{\mrd} f(t)\, dt=1,$ the
function $h$ has the partition of unity property with respect to the matrix $\cam.$
If the function $g$ is nonnegative, then the transform
$K_{g}$ maps nonnegative functions $f$
to nonnegative functions $h=K_g f.$\epr

\bp The assumptions imply that $g$ is bounded, so it is
clear that the integral in \eqref{51224bda} is well-defined for every $\ga \in \mrd$
whenever $f\in L^1(\mrd).$ Fixing any $\ga \in \mrd \setminus \{0\},$
\bes \inrd  \sum_{j\in \mz} \big| f(t) g\big(\frac{\cam^j \ga}{||t||} \big)\big|\,dt
=\inrd   | f(t)| \sum_{j\in \mz} |g(\cam^j( \ga ||t||^{-1} ))|\,dt
\le C\inrd |f(t)|\, dt
< \infty;\ens
thus, by Lebesgue's dominated convergence theorem,
\bes \sum_{j\in \mz} h(\cam^j \ga)=\inrd f(t) \sum_{j\in \mz}g(\cam^j( \ga ||t||^{-1} ))\,dt= \inrd f(t)\, dt.
\ens The rest of the proof is clear.\ep


A similar but more general result can be obtained by replacing the expression
$g(\frac{\ga}{||t||})$ in \eqref{51224bda} by a function $g(t,\ga)$ that yields
a partition of unity in the second variable. We leave the exact formulation to the
interested reader.

\subsection{An example of the integral transform on $\mr$ and a class of splines} \label{60129b}
In this subsection, we will study the one-dimensional version of the integral transform
in \eqref{51224bda}.
We will fix a constant $c\in (0,1),$ and consider the set
\bes 
S:= [-1, - c) \cup (c,1].\ens
Furthermore, we will fix $g:= \chi_S.$
Then the integral transform $K_g$
in \eqref{51224bda}, which we denote simply as $K$ here,  takes the form
\bee \label{60109d} h(\ga)= Kf(\ga):= \inr f(t) \chi_{S}\big(\frac{\ga}{|t|}\big)\,dt, \, \ga \in \mr. \ene
Note that for any fixed $\ga \in \mr,$
\bes 
\chi_S(\frac{\ga}{|t|})=1 \Leftrightarrow c\,|t| < | \ga| \le |t|
\Leftrightarrow  |\ga| \le |t| < |\ga| / c;\ens thus,
\begin{equation} \label{60109df} h(\ga)= \int_{-|\ga|/c}^{-|\ga|}f(t)\,dt+\int_{|\ga|}^{|\ga|/c} f(t)\,dt.\end{equation}
In particular,
the integral in \eqref{60109df}
is well-defined for all $\ga \in \mr$ whenever $f\in L_{\mbox{\tiny loc}}^1(\mr).$
We leave the short proof of the following result to the reader.

\bl \label{60115a} For any  $f\in L_{\mbox{{\rm {\tiny loc}}}}^1(\mr),$ the function $h=Kf$ is even; and if $f$
is an even function, then for $\ga >0,$
\bes
h(\ga) = 2 \int_{\ga}^{\ga/c} f(t)\, dt.\ens \el

The main merits of the transform $K$ are that it increases the regularity of $f$ and that the resulting function
$h=Kf$ satisfies the scaling partition of unity
property under some weak conditions on $f$:

\bpr \label{60109e} Let $f\in L_{\mbox{{\rm {\tiny loc}}}}^1(\mr) $ and consider the integral transform
$h=Kf$ in \eqref{60109d}. Then the following hold:
\bei \item[{\rm (i)}] If $f\in L^1(\mr),$ then
\bes
\sum_{j\in \mz} h(c^j\ga)= \inr f(t)\, dt, \, \forall \ga \in \mr \setminus \{0\}.\ens
\item[{\rm (ii)}] If $f\in C^k(\mr)$ for some $k\in \mn \cup \{0\}$ and  $f$ is supported away
from the origin, then $ h\in C^{k+1}(\mr).$ \eni\epr

\bp  As (i) clearly follows from Proposition \ref{51224eba},
we only have to prove (ii).  Letting
$F(\ga):= \int_0^\ga f(t)\, dt, \, \ga \in \mr,$
it follows from \eqref{60109df} that
\bes h(\ga) & = &
F( |\ga|/c)-F(|\ga|) + F(-|\ga|)-F(- |\ga|/c) \\ & = & \begin{cases}
F( \ga/c)-F(\ga) + F(-\ga)-F(- \ga/c), &\mbox{if} \ \ga \ge 0, \\
-[F( \ga/c)-F(\ga) + F(-\ga)-F(- \ga/c)], &\mbox{if} \ \ga \le 0.\end{cases}
\ens
For $\ga > 0,$ the function $h$ is obviously differentiable, and
\bes h^\prime (\ga)= \frac1{c} f(\ga/c)-f(\ga)-f(-\ga)+ \frac1{c} f(-\ga/c); \ens
thus under the stated assumptions $h$ is $(k+1)$ times continuously differentiable for $\ga >0.$ Similarly, $h$ is $(k+1)$
times continuously differentiable for $\ga<0;$ and since the function $h$ vanishes on a neighborhood of zero, $h$ is even infinitely differentiable at $\ga =0.$ \ep

\bex Let $f(t)=e^{-|t|}, \, t \in \mr.$ Then for $\ga \in \mr,$
\bes h(\ga)= \inr e^{-|t|} \chi_S\big(\frac{\ga}{|t|}\big)\,dt  =  2 \int_{|\ga|}^{|\ga|/c}e^{-|t|} \, dt
=  2( e^{-|\ga|} - e^{-|\ga|/c}).\ens
Observe that Proposition \ref{60109e}(i) implies that $h\in C(\mr)$ and
$\sum_{j\in \mz}h(c^j\ga)=2$ for $\ga \in \mr \setminus \{ 0 \}.$
We could of course obtain this construction via Proposition \ref{60108f}
as well.
\ep \enx

We will now use the integral transform $K$ to give a
direct definition of  a class of splines with attractive properties.

\bd \label{60109h} Let $h_1:= \chi_S,$ and define the functions $h_n, n\ge 2,$ inductively by
\bee \label{60109i} h_n(\ga):= Kh_{n-1}(\ga)=\inr h_{n-1}(t) \chi_{S}\big(\frac{\ga}{|t|}\big)\,dt, \, \ga \in \mr. \ene \ed

\bex \label{60115d} Direct calculation based on \eqref{60109i} shows that
\bes h_2(\ga)= \begin{cases} 0, &\mbox{if} \  |\ga| \le c^2, \\
2c^{-1} |\ga|-2c,  &\mbox{if} \  c^2 \le |\ga| \le c, \\
 2- 2\,|\ga|, &\mbox{if} \  c \le |\ga| \le 1, \\
 0, &\mbox{if} \ 1 \le  |\ga|,
  \end{cases} \ens
and
\bes h_3(\ga)= \begin{cases} 0, &\mbox{if} \  |\ga| \le c^3, \\
2c^{-3} |\ga|^2 - 4\, |\ga| +2c^3
, &\mbox{if} \  c^3 \le |\ga| \le c^2, \\
-2(c^{-1}+c^{-2})\,|\ga|^2 + 4(c+c^{-1})\, |\ga|- 2(c+c^2),  &\mbox{if} \  c^2 \le |\ga| \le c, \\
 2(1-|\ga|)^2, &\mbox{if} \  c \le |\ga| \le 1, \\
 0, &\mbox{if} \ 1 \le  |\ga|.
  \end{cases} \ens

\ep \enx

Let us collect some of the key properties of the spline functions $h_n$:

\bpr \label{60114a} The functions $h_n, \, n\in \mn,$ have the following properties:
\bei \item[{\rm (i)}] $h_n$ is a spline, with knots at the points $\pm c^n, \pm c^{n-1}, \dots, \pm 1.$
\item[{\rm (ii)}] $h_n$ is even.
\item[{\rm (iii)}] For $n \geq 2,$ $h_n\in C^{n-2}(\mr).$
\item[{\rm (iv)}] \, $\mbox{\rm supp} \, h_n= [-1, -c^n] \cup [c^n,1]$ and $h_n>0$ on
$(-1, -c^n) \cup (c^n,1).$
\item[{\rm (v)}] $Q_n:= \inr h_n(\ga)\, d\ga>0$ for all $n\in \mn,$ and $Q_1=2(1-c).$
\item[{\rm (vi)}] $ \frac{1}{Q_{n-1}}\, h_n$ satisfies the partition of unity condition
\bes 
\frac{1}{Q_{n-1}}\sum_{j\in \mz} h_n (c^j\ga)=1, \, \forall \ga \in \mr \setminus \{0\}.\ens
\item[{\rm (vii)}] There exists a constant $C>0$ such that
\bes 
C \le \frac{1}{Q_{n-1}^2}\sum_{j\in \mz} |h_n (c^j\ga)|^2\le 1, \, \forall \ga \in \mr \setminus \{0\}.\ens
\item[{\rm (viii)}] For $n\ge 2,$  the functions $h_n$ satisfy the recursion formula
\bee \label{60109j} h_n(\ga)= \frac{2}{n-1} \left[ (1-|\ga|) h_{n-1}(\ga)+
(c^{-1}|\ga| - c^{n-1})h_{n-1}(c^{-1}\ga)\right], \, \ga \in \mr.\ene
\eni
\epr

\bp Most of the results are immediate consequences of results that are already
proved. Indeed, (i) follows from (viii), which will be proved below; (ii) follows from the
definition and Lemma \ref{60115a}; and (iii) and (vi) are obtained from Proposition \ref{60109e}
and Example \ref{60115d}. In addition,
(iv) is proved by a straightforward induction, (v) is a
consequence of (iv) plus a direct calculation of $Q_1;$ and (vii) follows from
the partition of unity exactly as in the proof of Proposition \ref{60113d}(iii).

We will now prove the only item that remains, namely (viii). Since $h_n$ is even for
all $n\in \mn,$ we will assume that $\ga \ge 0.$
To get started, direct calculations based on the expressions in Example \ref{60115d}
show that the recursion formula holds for $n=2$ and $n=3.$  Thus, we will now consider $n\ge 4.$ Define the function $H_n$ by
\bee \label{60115fg} H_n(\ga):= \int_0^\ga h_n(t)\, dt, \, \ga \geq 0.\ene
We will perform an inductive proof of the recursion formula for $h_n,$ assuming
that it holds for $h_k$ for all $k=2, \dots, n-1.$ Now, using
Lemma \ref{60115a} and the induction hypothesis,
\bes
h_n(\ga) =
\frac{4}{n-2} \int_\ga^{\ga/c} \left[(1-t) h_{n-2}(t)+
(c^{-1}t - c^{n-2})h_{n-2}(c^{-1}t)    \right] \,dt.
\ens Then a direct calculation using integration by parts yields that
\bee \label{3.7aa} h_n(\ga) & = & \frac{4}{n-2} \bigg[ -\int_\ga^{\ga/c}(H_{n-2}(t/c)-H_{n-2}(t))\,dt+
(1-\ga)(H_{n-2}(\ga/c)-H_{n-2}(\ga)) \notag \\ 
& \ & \hspace{4cm} + \, (c^{-1}\ga-
c^{n-1})(H_{n-2}(\ga/c^2)-H_{n-2}(\ga/c))\bigg].\ene
Now, it follows from \eqref{60115fg} and Lemma \ref{60115a} that
\bes \int_\ga^{\ga/c}(H_{n-2}(t/c)-H_{n-2}(t))\,dt =
\frac12 \int_\ga^{\ga/c} h_{n-1}(t)\, dt = \frac14 h_n(\ga). \ens
Also,
$$ H_{n-2}(\ga/c)-H_{n-2}(\ga)= \int_\ga^{\ga/c} h_{n-2}(t)\, dt = \frac12 h_{n-1}(\ga).$$
Hence, based on \eqref{3.7aa}, after solving for $h_n(\ga)$, we obtain \eqref{60109j}.
\ep

The splines in Definition \ref{60109h} are
indeed well-known: as noted from the recursion formula \eqref{60109j}, they are the symmetrized version of the
nonuniform B-splines with
knots at $c^n, c^{n-1}, \dots, 1$, see \cite{Boor, BHR}. Here, we have provided another perspective in obtaining them.
Their properties also serve as a concrete illustration of the general properties we derived in
Propositions \ref{51224eba} and \ref{60109e}.
Other related papers on polynomial splines with geometric knots include \cite{GLee1, Lee1, Mic}.


As a further comment on the one-dimensional transform $K$ in \eqref{60109d}, we observe that it can be lifted to a transform acting on functions
on $\mrd$:

\bex \label{60214b}
In this example, we describe a way of lifting the transform $K$ to generate radial functions
on $\mrd.$

\vspace{.1in}\noindent(i)
We can easily lift the integral transform to an operator that yields a radial function
$\widetilde{h}: \mrd \to \mr$ as output.
Indeed, taking an arbitrary norm $|| \cdot ||$ on $\mrd,$ define the integral transform $\widetilde{K},$ acting on functions $f\in L_{\mbox{\tiny loc}}^1(\mr)$, by
\bes 
\widetilde{h}(\ga)= \widetilde{K}f(\ga):= \int_{-\infty}^\infty f(t) \chi_{S}\big(\frac{||\ga||}{|t|}\big)\,dt, \, \ga \in \mrd. \ens
Clearly, in terms of the transform $K$  in \eqref{60109d},  we have $\widetilde{h}(\ga)=Kf(||\ga||) = h(|| \ga ||).$ Furthermore, if $f\in L^1(\mr),$ then
\bes  \sum_{j\in \mz} \widetilde{h}(c^j\ga)= \inr f(t)\, dt, \, \forall \ga \in \mrd \setminus \{0\}.\ens

\vspace{.1in}\noindent(ii) As a special case of (i)
and based on the nonuniform B-splines $h_n$ in \eqref{60109i},  we can define a family of radial functions $\widetilde{h_n}$ on $\mrd$ by
\bes \widetilde{h_n}(\ga):=h_n(||\ga||), \, \ga \in \mrd. \ens
Each of these radial functions is supported on an annulus, and they can be easily
calculated using the recursion formula in Proposition \ref{60114a}. Also, $\widetilde{h_n}$
satisfies the partition of unity condition
$$
\sum_{j\in \mz} \widetilde{h_n}(c^j\ga)= \int_{-\infty}^\infty h_{n-1}(t) \, dt, \, \forall \gamma \in \mrd \setminus \{0 \}.
$$
\ep \enx


\section{Wavelet frames in $L^2(\mrd)$ and dual frames} \label{60129c}
In this section, we apply the results on the scaling partition of unity to
construct dual pairs of matrix-based wavelet frames in $\ltrd.$ Since wavelet frames
is a well-studied area by itself (see, e.g., \cite{RoSh2,DRoSh3,CHS}), we will not make any attempt to motivate them or highlight
their applications but just state the definitions and results that are strictly necessary
for our discussion. 

Given an invertible $d\times d $ matrix $\cam$ with real entries,
we define the scaling operator $D_\cam: \ltrd \to \ltrd$ by
$ (D_\cam f)(x):=| \det \cam |^{1/2}f(\cam x);$ and, for  $\nu\in \mrd,$ the
translation operator $T_\nu: \ltrd\to \ltrd$ by $T_\nu f(x):= f(x-\nu).$
Fixing a function $\psi\in \ltrd,$ a $d\times d$  matrix $\cam$
and a translation parameter $b>0,$ the associated wavelet system is given by
$\{ D_{\cam^j} T_{bk}\psi\}_{j\in \mz, k\in \mzd}.$ Denoting the
canonical norm on $\ltrd$ by $||\cdot ||_2,$ the wavelet system $\{ D_{\cam^j} T_{bk}\psi\}_{j\in \mz, k\in \mzd}$ is said to form a frame
for $\ltrd$ if there exist constants $A,B>0$ such that
\bee \label{60121c} A ||f||_2^2 \le \sum_{j\in \mz, k\in \mzd} | \la f, D_{\cam^j} T_{bk}\psi \ra |^2 \le B ||f||_2^2, \, \forall f\in \ltrd;\ene if at least the upper condition in \eqref{60121c} is satisfied, it is called
a Bessel sequence. Two Bessel sequences $\{ D_{\cam^j} T_{bk}\psi\}_{j\in \mz, k\in \mzd}$ and $\{ D_{\cam^j} T_{bk}\widetilde{\psi}\}_{j\in \mz, k\in \mzd},$
where $\psi, \widetilde{\psi} \in \ltrd,$ are said to form
dual frames if
\bes 
f= \sum_{j\in \mz, k\in \mzd}  \la f, D_{\cam^j} T_{bk}\widetilde{\psi} \ra D_{\cam^j} T_{bk}\psi, \, \forall f\in \ltrd.\ens

We will need the
following result, which gives sufficient conditions for wavelet systems to form
Bessel sequences, frames, and dual frames.
It exists in several variants in the literature: (i) was first
stated explicitly in \cite{Lem2}, while versions of (ii) can be found, e.g., in
\cite{ChuiCzaja,HLW}; see also \cite{CBN}. We define the Fourier transform on
$L^1(\mrd)$ by ${\cal F}f(\ga)= \widehat{f}(\ga):= \int_{\mr^d} f(x) e^{-2\pi i x \cdot \ga}\, dx,
\, \ga \in \mrd,$
with the usual extension to $\ltrd.$

\bl \label{41131a} Let $\cam$ denote an invertible $d\times d$ matrix with real entries, and
let $b>0.$ Then the following hold:

\bei \item[{\rm (i)}] If $\psi\in \ltrd$ and
$B:= \frac1{b^d}\, \esssup_{\ga \in \mrd} \sum_{j\in \mz} \sum_{k\in \mzd}
| \widehat{\psi}((\cam^T)^j \ga) \, \widehat{\psi}((\cam^T)^j \ga
- k/b)|< \infty,$ then $\{D_{\cam^j}T_{bk}\psi\}_{j\in \mz, k\in \mzd}$ is a Bessel sequence.
If furthermore
\bes A:= \frac1{b^d}\, \essinf_{\ga \in \mrd} \left(
\sum_{j\in \mz}| \widehat{\psi}((\cam^T)^j\ga)|^2 -
\sum_{j\in \mz} \sum_{k\neq 0}
| \widehat{\psi}((\cam^T)^j \ga) \, \widehat{\psi}((\cam^T)^j \ga
- k/b)|\right)>0,\ens then $\{D_{\cam^j}T_{ bk}\psi\}_{j\in \mz, k\in \mzd}$ is a frame for $\ltrd$ with
bounds $A, B.$
\item[{\rm (ii)}] Assume that the matrix $\cam$ is expanding and suppose that for some $\psi, \widetilde \psi \in \ltrd,$
$\{D_{\cam^j}T_{ bk}\psi\}_{j\in \mz, k\in \mzd},$ $\{D_{\cam^j}T_{ bk}\widetilde{\psi}\}_{j\in \mz, k\in \mzd}$
are Bessel sequences.
Then $\{D_{\cam^j}T_{ bk}\psi\}_{j\in \mz, k\in \mzd},$ $\{D_{\cam^j}T_{ bk}\widetilde{\psi}\}_{j\in \mz, k\in \mzd}  $ are
dual frames for $\ltrd$ if and only if for all $m \in \mzd,$
\bee \label{60120a} \sum_{\{ j \in \mz \,  | \,   (A^T)^{-j} m\in \mzd   \}}
\overline{\widehat{\psi}((\cam^T)^{-j} \ga)}
\, \widehat{\widetilde{\psi}}((\cam^T)^{-j} \ga + (\cam^T)^{-j} m/b)
=b^d \, \delta_{m,0}, \, \mbox{a.e.}\ \ga \in \mrd.\ene \eni \el

\begin{rem} \label{60129d}
\rm{It follows from Lemma \ref{41131a}  that if $\widehat{\psi}$ is supported
on the closed ball $B(0,R)$ of radius $R$ in $\mrd$ and $b\le (2R)^{-1},$ then $\{D_{\cam^j}T_{bk}\psi\}_{j\in \mz, k\in \mzd}$ is a Bessel sequence when $\esssup_{\ga \in \mrd} \sum_{j\in \mz}
| \widehat{\psi}((\cam^T)^j \ga)|^2< \infty;$ and it is a frame when
\bes 
0 <\essinf_{\ga \in \mrd} \sum_{j\in \mz}
| \widehat{\psi}((\cam^T)^j \ga)|^2\le
\esssup_{\ga \in \mrd} \sum_{j\in \mz}
| \widehat{\psi}((\cam^T)^j \ga)|^2< \infty.\ens
If both $\widehat{\psi}$ and  $\widehat{\widetilde{\psi}}$ are supported
on  $B(0,R)$, then
\eqref{60120a} is satisfied for
$m \in \mzd \setminus \{0\}$ when  $b \leq (2R)^{-1};$ in this case the
condition \eqref{60120a} consists of the single equation
\bes 
\sum_{ j \in \mz}
\overline{\widehat{\psi}((\cam^T)^{-j} \ga)}
\, \widehat{\widetilde{\psi}}((\cam^T)^{-j} \ga)
=b^d, \, \mbox{a.e.}\, \ga \in \mrd.\ens Up to the factor $b^d,$ this essentially means that the function $\overline{\widehat{\psi}}\widehat{\widetilde{\psi}}$ satisfies the scaling partition
of unity property with respect to the matrix $\cam^T.$}
\end{rem}

Proposition \ref{60114a} and Lemma \ref{41131a} lead to the following frame result on $\ltr$
for the splines $h_n$ in \eqref{60109i}:

\bt \label{60129g} Given any $n\in \mn$ and $c\in (0,1),$ consider the spline $h_n$ in  \eqref{60109i}. Fix $b\in (0, c^{n-1}/2]$
and define the functions $\psi, \widetilde{\psi} \in \ltr$ by $\widehat{\psi}:= h_n$ and
\bee \label{60121f} \widehat{ \widetilde{\psi}} (\ga):= \frac{b}{Q_{n-1}^2}\sum_{j=-n+1}^{n-1}h_n (c^j\ga), \, \ga \in \mr.\ene Then $\{D_{c^j}T_{kb}\psi\}_{j,k\in \mz}$ and $\{D_{c^j}T_{kb}\widetilde{\psi}\}_{j,k\in \mz}$
are dual wavelet frames for $\ltr.$
\et

\bp Since $\supp\, \widehat{\psi} \subset [-1,1] \subseteq [-c^{-n+1}, c^{-n+1}],$ the frame property of
$\{D_{c^j}T_{kb}\psi\}_{j,k\in \mz}$
follows directly from
Proposition \ref{60114a}(vii) and Remark \ref{60129d}. Now, by the partition of unity condition in
Proposition \ref{60114a}(vi), we have
\bee \label{60121g}  \frac{1}{Q_{n-1}}\sum_{j\in \mz} \widehat{\psi} (c^j\ga)=1, \, \forall \ga \in \mr \setminus \{0\}.\ene The expression on the right-hand side of \eqref{60121f} clearly
defines a bounded function, with compact support
$[-c^{-n+1}, -c^{n-1}] \cup [c^{n-1}, c^{-n+1}]$ which is bounded away from the origin. Thus the
function $\widetilde{\psi}$ is well-defined and $\{D_{c^j}T_{kb}\widetilde{\psi}\}_{j,k\in \mz}$ is a Bessel sequence by Lemma \ref{41131a}(i) and Remark \ref{60129d}.

If $\ga \in \supp\, \widehat{\psi}$, then $\widehat{\psi}(c^j\ga)$
can only be nonzero for $j=-n+1, -n+2, \dots, n-1;$ thus \eqref{60121g}
implies that $\widehat{\widetilde{\psi}}(\ga)=\frac{b}{Q_{n-1}}$ for $\ga \in \supp\, \widehat{\psi}.$
It follows that $ \overline\,
\widehat{\widetilde{\psi}}(\ga)   = \frac{b}{Q_{n-1}} \overline{\widehat{\psi} (\ga)}$ for all $\ga \in \mr;$
using again \eqref{60121g} now
shows that
\bes  \sum_{j\in \mz} \overline{\widehat{\psi} (c^j\ga)}
\, \widehat{\widetilde{\psi}}(c^j\ga)   =b, \, \forall \ga \in \mr \setminus \{0\}.\ens
Hence we conclude from Lemma \ref{41131a}(ii) and Remark \ref{60129d} that
$\{D_{c^j}T_{kb}\psi\}_{j,k\in \mz},$
$\{D_{c^j}T_{kb}\widetilde{\psi}\}_{j,k\in \mz}$ are indeed dual frames.\ep

Note that a different dual frame $\{D_{c^j}T_{kb}\widetilde{\psi}\}_{j,k\in \mz}$
associated with $\{D_{c^j}T_{kb}\psi\}_{j,k\in \mz}$ could have been obtained via
the results in \cite{Lem1}. Also, by combining the result with the lifting
transform in Example \ref{60214b}, it is easy to construct radial dual wavelet
frames $\{D_{c^jI}T_{kb}\psi\}_{j\in \mz, k \in \mzd}$
and $\{D_{c^jI}T_{kb}\widetilde{\psi}\}_{j\in \mz, k \in \mzd}$ for $\ltrd,$ where $\psi, \widetilde{\psi} \in \ltrd;$ we leave the
details to the reader.

We also note that the unitary extension principle 
and its many variants is a classical tool to construct wavelet frames based on splines,
see, e.g., \cite{RoSh2,DRoSh3,CHS}.
However, in this case the frame generators themselves are splines, while in our
construction the splines occur in the Fourier domain.

We will now establish a
result about the construction of wavelet frames in $\ltrd$,  based on Proposition \ref{60113d}:

\bt \label{60113da} Let $|| \cdot ||$ denote any norm on $\mrd,$ and consider an expanding $d \times d$
matrix $\cam$ such that $||\ga||\le ||\cam^T \ga||$ for all $\ga \in \mrd.$
Let $r: [0,\infty) \to \mr$ be a continuous decreasing function
supported on $[0,R]$ for some $R>0$ such that
$r(0)=1$.  Consider the function $\psi: \mrd \to \mr$
defined via
$\widehat{\psi}(\ga):=r(||\ga||)-r(|| \cam^T\ga||), \, \ga \in \mrd.$
Then the following hold:
\bei \item[{\rm (i)}] Whenever $b\le (2R)^{-1},$ $\{D_{\cam^j}T_{ bk}\psi\}_{j\in \mz, k\in \mzd}$
is a wavelet frame for $\ltrd$.
\item[{\rm (ii)}] If $r(\ga)=1$ for $\ga\in [0, R_1]$ for some $R_1>0,$
there exists a finite index set $J$ containing $0,$ which depends on the
matrix $\cam$ and the numbers $R,R_1,$ such that
\bes 
\sum_{j\in J} \widehat{\psi}((\cam^T)^j \ga)=1, \, \forall
\ga \in \mbox{\rm supp} \, \widehat{\psi}.\ens
\item[{\rm (iii)}]
If $r(\ga)=1$ for $\ga\in [0, R_1]$ for some $R_1>0,$ choose an index set
$J$ as in {\rm (ii)} and for $b > 0,$
define the function $\widetilde{\psi}: \mrd \to \mr$ via
\bee \label{4.4aa}  \widehat{ \widetilde{\psi}} (\ga):= b^d \sum_{j\in J} \widehat{\psi}((\cam^T)^j \ga), \, \ga \in \mrd.\ene
Then for sufficiently small values of $b,$
 $\{D_{\cam^j}T_{ bk}\psi\}_{j\in \mz, k\in \mzd}$ and $\{D_{\cam^j}T_{ bk}\widetilde{\psi}\}_{j\in \mz, k\in \mzd}  $ are
dual frames for $\ltrd$. \eni \et

\bp The matrix $\cam^T$ is expanding, so Proposition \ref{60113d}(iii) implies that there
exists a constant $C>0$ such that
$C \le \sum_{j\in \mz} | \widehat{\psi}((\cam^T)^j\ga)|^2 \le 1$ for all $\ga \in \mrd \setminus \{0\}.$ The result in (i) now follows from Remark \ref{60129d}.

In order to prove (ii), we will first show that the function $\widehat{\psi}$ is supported on
the annulus $a(R_1 || \cam^T||^{-1}, R).$ (This annulus is well-defined as $|| \ga || \leq || \cam^T \ga || \leq ||\cam^T || \,
||\ga||$ implies that $|| \cam^T || \geq 1.$) If $||\ga||\ge R,$ we also have that
$||\cam^T\ga||\ge R,$ so indeed $\widehat{\psi}(\ga)=0.$ Now, assume that
$||\ga|| \le R_1 || \cam^T||^{-1}.$ Then $||\ga|| \le R_1$ and
\bes ||\cam^T \ga|| \le||\cam^T|| \, || \ga|| \le ||\cam^T|| \, R_1 || \cam^T||^{-1} =  R_1; \ens
thus $\widehat{\psi}(\ga)= r(||\ga||)-r(|| \cam^T\ga||)=1-1=0$
as claimed.

Now, applying Proposition \ref{60113d}(ii) to the expanding
matrix $\cam^T,$ we have
\bee \label{60129h} \sum_{j\in \mz} \widehat{\psi}((\cam^T)^j\ga) = 1, \, \forall \ga \in \mrd \setminus \{0\}. \ene
Since $\cam^T$ is expanding, there
exist constants $C' \in (0,1]$ and $\alpha>1$ such that
$|| (\cam^T)^j \ga || \ge C' \alpha^j ||\ga||$ for all $\ga\in \mrd$ and
$j\in \mn \cup \{0\}.$ Thus, for $j\in \mn \cup \{0\}$ and any $\ga\in \supp\, \widehat{\psi},$
\bes || (\cam^T)^j \ga || \ge C' \alpha^j ||\ga|| \ge  C' \alpha^j R_1 || \cam^T||^{-1}. \ens
It follows that for $\ga \in \supp\, \widehat{\psi},$
we have $\widehat{\psi}((\cam^T)^j \ga)=0$ whenever $C' \alpha^j R_1 || \cam^T||^{-1}
\ge R,$ i.e., for $j$ sufficiently large.  On the other hand, since $||(\cam^T)^{-j}\ga|| \le (C')^{-1} \alpha^{-j}||\ga||$
for all $\ga\in \mrd$ and $j\in \mn \cup \{0\},$ it also follows that for
$\ga \in \supp \, \widehat{\psi},$
\bes ||(\cam^T)^{-j}\ga|| \le (C')^{-1} \alpha^{-j}R;\ens thus
$\widehat{\psi}((\cam^T)^{-j} \ga)=0$ whenever $(C')^{-1} \alpha^{-j}R \le R_1 || \cam^T||^{-1},$
i.e., for $j$ sufficiently large. This completes the proof of (ii).

Finally, to establish (iii), we first note that for every $j \in \mz,$ there exist constants $\lambda_j, \mu_j > 0$ such that
\bee
\label{4.6a}
 \lambda_j || \ga || \leq || (\cam^T)^j \ga || \leq \mu_j || \ga ||, \, \forall \ga \in \mrd.
\ene
Indeed, for $j=0,$ we simply take $\lambda_ 0 = \mu_0 =1.$
If $j \in \mn,$ then $C' \alpha^j ||\ga|| \leq || (\cam^T)^j \ga || \leq || (\cam^T)^j || \, || \ga ||$ for all $\ga\in \mrd.$
On the other hand, $|| (\cam^T)^{-j} \ga || \leq (C')^{-1} \alpha^{-j} ||\ga||$ and
$|| \ga || = || (\cam^T)^{j} (\cam^T)^{-j} \ga || \leq || (\cam^T)^{j} || \, || (\cam^T)^{-j} \ga ||$
for all $\ga\in \mrd.$ Next, using the fact that $\supp \, \widehat{\psi} \subset a( R_1 || \cam^T ||^{-1}, R),$
it follows from \eqref{4.6a} that $\supp \, \widehat{\psi}((\cam^T)^j \cdot) \subset a( R_1 || \cam^T ||^{-1} \mu_j^{-1}, R \lambda_j^{-1})$
for every $j \in \mz.$ Thus the definition of $\widehat{\widetilde{\psi}}$ in \eqref{4.4aa} shows that
$\widehat{\widetilde{\psi}}$ is a bounded function, with
\bes
\supp \, \widehat{\widetilde{\psi}} \subset \bigcup_{j \in J}
\supp \, \widehat{\psi}((\cam^T)^j \cdot) \subset
a\big( R_1 || \cam^T ||^{-1} \min_{j \in J} \mu_j^{-1}, R \max_{j \in J} \lambda_j^{-1} \big).
\ens
Since $J$ is a finite set containing $0,$ $\max_{j \in J} \lambda_j^{-1} \geq \lambda_0^{-1} = 1$ and so
$R \max_{j \in J} \lambda_j^{-1} \geq R.$
Consequently, both $\widehat{\psi}$ and $\widehat{\widetilde{\psi}}$ are supported on the closed ball
$B(0, R \max_{j \in J} \lambda_j^{-1} ).$ The rest of the proof of (iii) is similar to the proof of Theorem \ref{60129g},
where Lemma \ref{41131a} and Remark \ref{60129d} are applied. Specifically, we see that whenever
$b \leq (2R \max_{j \in J} \lambda_j^{-1})^{-1}$,
$\{D_{\cam^j}T_{ bk}\psi\}_{j\in \mz, k\in \mzd}$ and $\{D_{\cam^j}T_{ bk}\widetilde{\psi}\}_{j\in \mz, k\in \mzd}$
are Bessel sequences. Also, the partition of unity condition \eqref{60129h} together
with (ii) shows that
\bes \sum_{j\in \mz} \overline{\widehat{\psi}((\cam^T)^j\ga)} \,
\widehat{\widetilde{\psi}}((\cam^T)^j\ga)= b^d, \, \forall \ga \in \mrd \setminus \{0\};\ens
and hence, $\{D_{\cam^j}T_{ bk}\psi\}_{j\in \mz, k\in \mzd}$ and $\{D_{\cam^j}T_{ bk}\widetilde{\psi}\}_{j\in \mz, k\in \mzd}$
are dual frames for $\ltrd$. \ep


\medskip

\vspace{.1in}\noindent{\bf Acknowledgment:} Ole Christensen would like to thank the National University of Singapore for its warm hospitality during one-month stays in 2016 and 2017.

\begin{tabbing}
text-text-text-text-text-text-text-text-text-text \= text \kill \\
Ole Christensen \> Say Song Goh \\
Technical University of Denmark \> Department of Mathematics \\
DTU Compute \> National University of Singapore \\
Building 303, 2800 Lyngby \> 10 Kent Ridge Crescent \\
Denmark \> Singapore 119260, Republic of Singapore \\
Email: ochr@dtu.dk \> Email: matgohss@nus.edu.sg
\end{tabbing}

\end{document}